{

\documentclass[12pt]{article}
\usepackage[pdftex]{graphicx}
\usepackage{amsfonts}
\usepackage{mathrsfs}
\newcommand{\R}{\mathbb{R}}
\newcommand{\CC}{{\cal C}}
\newcommand{\Po}{\mathscr{P}}
\renewcommand{\epsilon}{\varepsilon}

\newtheorem{theo}{Theorem}

\newtheorem{prop}{Proposition}
\newtheorem{lemm}{Lemma}

\newcommand{\bea}{\begin{eqnarray}}
\newcommand{\eea}{\end{eqnarray}}
\newcommand{\beaa}{\begin{eqnarray*}}
\newcommand{\eeaa}{\end{eqnarray*}}

\title{Strict inequalities of critical probabilities on 
Gilbert's continuum percolation graph}

\author{Massimo Franceschetti\footnote{Electrical and Computer Engineering,
University of California San Diego,
La Jolla, CA, 92093-0407 
Email: massimo@ece.ucsd.edu},
Mathew D. Penrose\footnote{ Department of Mathematical Sciences,
 University of Bath BA2 7AY, United Kingdom,
Email: m.d.penrose@bath.ac.uk}
 and Tom Rosoman\footnote{ Department of Mathematical Sciences,
 University of Bath BA2 7AY, United Kingdom,
Email: ter20@bath.ac.uk}
}

\date{\today}

\begin{document}

\maketitle

\begin{abstract}\addcontentsline{toc}{section}{Abstract}
Any infinite graph has site and bond percolation critical probabilities   
satisfying 
 $p_c^{\rm site}\geq  p_c^{\rm bond}$. The strict version of this inequality holds for many, but not all, infinite graphs.

In this paper, the class of graphs for which the strict
inequality holds is extended to a continuum percolation model.
In Gilbert's graph with supercritical density on the Euclidian plane, there is almost surely a unique infinite connected component. We show that on this component $p_c^{\rm site} >  p_c^{\rm bond}$.
 This also holds in higher dimensions.
\end{abstract}

\section{Introduction}
\label{secintro}
Consider an infinite connected graph $G$ and perform bond
percolation by independently marking each edge open with probability
$p$ and closed otherwise. The critical probability $p_c^{\rm bond}$
refers to the value of $p$ above which there exists almost surely
(a.s.) an infinite connected subgraph of $G$, of open edges.
Similarly, one can perform site percolation by independently marking
each vertex of $G$ open with probability $p$ and refer to $p_c^{\rm
site}$ as the critical probability above which there exists a.s.\ an
infinite connected subgraph of $G$, of open vertices.

The weak inequality $p_c^{\rm site}\geq  p_c^{\rm bond}$ can easily be proven by dynamic coupling, see for example Chapter~2 of Franceschetti and Meester~(2007).  If $G$ is a tree,
% where the starting vertex is automatically open,
 then it is also easy to see that $p_c^{\rm site} =  p_c^{\rm bond}$, as each
vertex, other than some arbitrarily selected root, can be uniquely identified by an edge and vice versa. By adding finitely many edges to an infinite tree, one can also construct other connected graphs for which the equality holds.  

On the other hand, the strict 
%version of the
 inequality $p_c^{\rm site} >  p_c^{\rm bond}$ has also been shown to hold in several circumstances. Grimmett and Stacey~(1998) proved it for
a large class of `finitely transitive' graphs including the
  $d$-dimensional hypercubic lattices.

 These graphs, however, do not include the {\em random} graphs constructed 
using {\em continuum} percolation models,
%such as Gilbert's, 
%or the random connection model, 
because they are not `finitely transitive':
% in the sense of Grimmett and Stacey~(1998):
  since their average node degree is not bounded, the group action defined by their automorphisms has infinitely many orbits almost surely.  These continuum percolation graphs are the focus of this paper. They are of particular  interest in the context of communication networks and  are treated extensively in the books by  Franceschetti and Meester~(2007), Meester and Roy~(1996), and Penrose~(2003).

We consider {\em Gilbert's  graph},
 which is defined as follows. Let $\lambda >0$
and let $\mathscr{P}_\lambda$ be a homogeneous
 Poisson point process in $\R^2$ of intensity $\lambda$.
 Gilbert's graph, here denoted $G(\mathscr{P}_\lambda,1)$, 
has as its vertex set
 the point set $\mathscr{P}_\lambda$,
and the edges are obtained by connecting every 
 pair  of points $x,y \in \mathscr{P}_\lambda$
such that $|x-y|\leq 1$, 
by an undirected edge.
It is well known that there exists 
a critical density value $\lambda_c \in (0, \infty)$,
  such that if $\lambda > \lambda_c$ then there exists a.s.\ a unique infinite connected component, while  if $\lambda < \lambda_c$ then there is a.s no infinite connected component; 
see  e.g. Meester and Roy (1996).
 When it exists, we 
denote this infinite component by $\CC$.
 
In the site percolation model on  $\CC$, each vertex is independently marked open with probability $p$, and closed otherwise, and 
we look for an unbounded connected component in  the induced subgraph  $\CC_v$ of the open vertices.  It is easy to see that this is equivalent to rescaling the original Poisson process to one with intensity $p \lambda$ and looking for an unbounded connected component there. It follows that for $\lambda > \lambda_c$ there is a critical value $p_c^{\rm site} \in (0,1)$ 
(namely $p_c^{\rm site} = \lambda_c / \lambda$) such that if $p > p_c^{\rm site}$ then there is a.s.\ an infinite connected component in $\CC_v$, and if
 $p < p_c^{\rm site}$ then there is a.s.\ no such infinite component.

In the  bond percolation model on $\CC$, we independently declare each edge
to be open with probability $p$, and closed otherwise, and look for an
 unbounded connected component in the induced subgraph $\CC_e$ of the open
 edges.
% This is equivalent to constructing an RCM where the original
% connection function $f$ is replaced by $p \cdot f$.
There is a critical probability $p_c^{\rm bond} $ such that if $p > p_c^{\rm bond}$ then there is a.s.\ an infinite connected component in $\CC_e$, and
 if $p < p_c^{\rm bond}$ then there is 
a.s.\ no such  infinite component. 
Observe that $p_c^{\rm bond} \leq p_c^{\rm site} <1$, and it can also
be shown by a branching process comparison that 
$p_c^{\rm bond} >0  $.

Our main result provides strict inequality between 
$p_c^{\rm site}$ and $p_c^{\rm bond}$ on
 Gilbert's graph.  Our proof easily extends to $3$ or more dimensions.
\begin{theo} \label{bool}
Consider  $G({\mathscr P}_ \lambda,1)$ for $\lambda>\lambda_c$.
On $\CC$ we have $p_c^{\rm site} > p_c^{\rm bond}$.
\end{theo}

The basic strategy is to adapt the enhancement
 technique developed for percolation on lattices by
 Menshikov~(1987),
 Aizenman and Grimmett~(1991),
 Grimmett and Stacey~(1998).
 This consists of constructing an `enhanced' version of the site
 percolation process for which the critical probability is 
\emph{strictly less} than that of the original site process.
 Then one can use dynamic coupling of the enhanced model with
 bond percolation to complete the proof.
%  They also showed that a similar strategy, using a `diminishment' of bond percolation rather than an enhancement of site percolation, can be used to prove the strict inequality for a wide range of graphs beyond hypercubic lattices.

We face two main difficulties when trying to extend the enhancement technique
% results
 to a continuum random setting. One of these  amounts to constructing
 the desired enhancement on a random graph rather than on a deterministic
 one. The second one consists in adapting some basic inequalities for
 the enhanced graph, given in the discrete setting by Aizenman and
 Grimmett~(1991), to the continuum setting. This requires somehow
 more involved geometric constructions and a careful incremental 
build-up of the Poisson point process. Once we circumvent these obstacles,
 it is not too difficult to obtain the final result using a classic dynamic 
coupling construction.

The enhancement strategy has been proven useful to show strict inequalities in a variety of contexts: Bezuidenhout, Grimmett, and Kesten~(1993), and Grimmett~(1994), use this technique in the context of Potts and random cluster models;  Roy, Sarkar, and White~(1998) use it in the context of directed percolation.
In the continuum, Sarkar~(1997) uses enhancement to demonstrate coexistence
of occupied and vacant phases for the three-dimensional
  Poisson Boolean model.
  Roy and Tanemura~(2002) use it in the context of percolation of different convex shapes.

\section{Proof of Theorem \ref{bool}}
\begin{figure}[htbp]
\begin{center}
\scalebox{.8}{\includegraphics[angle = 270]{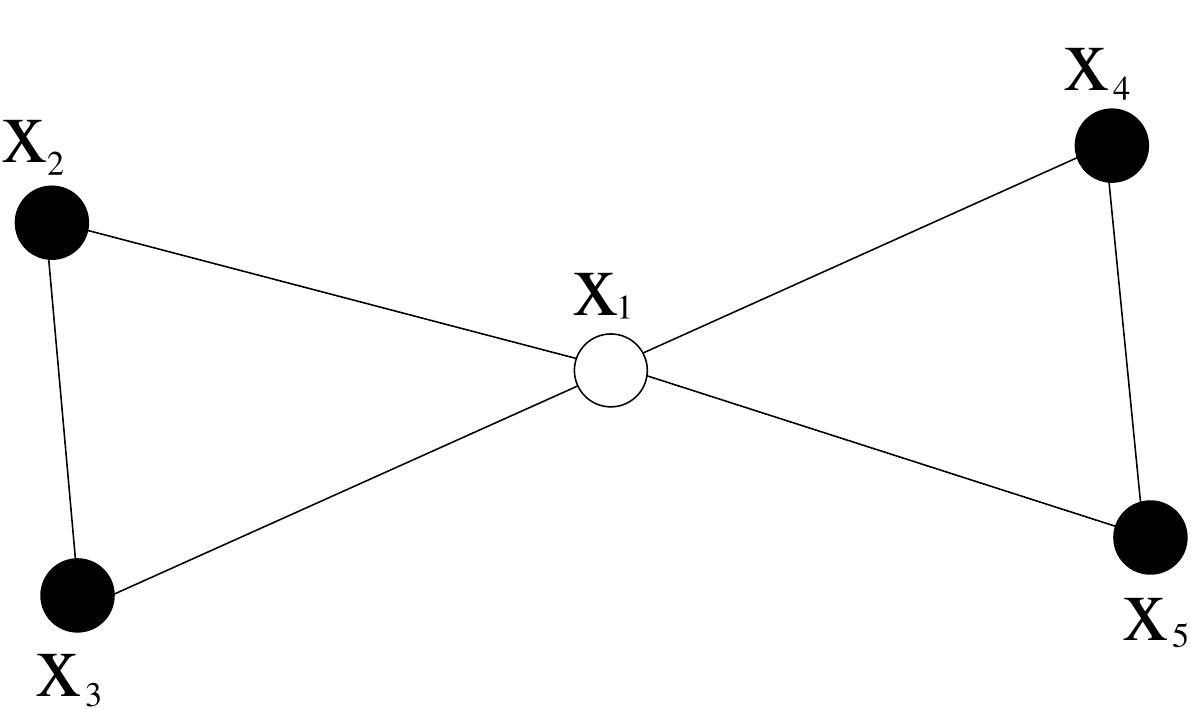}}
\end{center}
\caption{The bow tie enhancement}
\label{fig:enh1}
\end{figure}

We now describe the enhancement needed to prove Theorem $1$.
 Throughout this section we consider Gilbert's graph 
$G({\mathscr P}_ \lambda,1)$ with $\lambda > \lambda_c$.
The
objective is to describe a way to to add open vertices to the site
percolation model to make the probability of an infinite cluster
bigger, without changing the bond percolation model. To do so, we
introduce two kinds of coloured vertices, red vertices (the original
open vertices) and green vertices (closed vertices which have been
enhanced) and for any two vertices $x,y$ we write that $x \sim y$ if
they are joined by an edge. In $G(\mathscr{P}_\lambda,1)$, if we have vertices
$x_1,x_2,x_3,x_4,x_5$ such that $x_1$ is closed, has no neighbours
other than $x_2,\dots,x_5$, which are all red, and $x_2\sim x_3$ and
$x_4\sim x_5$ but there are no other edges amongst $x_2, x_3, x_4$
and $x_5$ then we say $x_1$ is {\em correctly configured}
 in $G(\mathscr{P}_\lambda,1)$, and refer to this as a $bow$ $tie$ configuration of edges. If a vertex $x$ is correctly configured
we make it green with probability $q$, independently of everything
else; see Figure~\ref{fig:enh1}.

Let $B_n$ be the open disc of radius $n$ centred at the origin. Let
$\underline{Y} = (Y_i, i \geq 0)$  and $\underline{Z} = (Z_i, i \geq
0)$ be sequences of independent uniform $[0,1]$ random variables.
List the vertices of $\mathscr{P}_{\lambda}$ in order of increasing distance
from the origin as $x_1,x_2,x_3,\dots$ . Declare a vertex $x_i$ to be 
{\em red} if
$Y_i < p$ and {\em closed} otherwise. Once the sets of red and closed
vertices have been decided in this way, apply the enhancement by
declaring each closed vertex $x_j$ to be {\em green}
% (and therefore coloured)
 if it is
correctly configured and $Z_j < q$.  Whenever we
insert a vertex of the Poisson process at $x$, it would have values $Y_0$ and $Z_0$ associated with it. We shall refer to vertices that are either
red or green as being {\em coloured}.

Let $\partial B_n$ be the annulus $B_n \setminus B_{n - 0.5}$ and let
$A_n$ be the event that there is a path from a coloured vertex in
$B_{0.5}$ to a coloured vertex in $\partial B_n$ in $G
(\mathscr{P_{\lambda}},1)\cap B_n$  using only coloured
vertices (note that $A_n$ is based on a process completely
 inside $B_n$; we do not allow vertices outside of $B_n$ to
 affect possible enhancements inside $B_n$). 

Let $\theta_n(p,q)$ be the probability that $A_n$ occurs, and define
$$\theta(p,q) \equiv \liminf _{n\rightarrow \infty} (\theta_n(p,q)).$$

The following proposition states that $\theta(p,q)$ is indeed the
percolation function associated to the enhanced model. From now on
we use `vertex' to refer to a point of the Poisson process and `point' to
refer to an arbitrary location in $\mathbb{R}^2$.

\begin{prop}\label{prop:enh}
There is a.s. an infinite connected component in $G({\mathscr P}_ \lambda,1)$ using only red and green vertices if and only if $\theta(p,q) > 0$.
\end{prop}
\noindent {\bf Proof of Proposition~\ref{prop:enh}.} For the if part let $A'_n$ be the event that there is a coloured path from $B_{0.5}$ to outside $B_{n-2}$, so $A_n$ is contained in $A'_n$. Let $\phi_n(p,q)$ be the probability of $A'_n$ occurring and let $\phi(p,q)$ be the limit as $n$ goes to $\infty$. Therefore $\phi_n(p,q) \geq \theta_n(p,q)$ for all $n$ so $\phi(p,q) \geq \theta(p,q) > 0$, but $\phi(p,q)$ is just the probability of there being an infinite coloured component intersecting $B_{0.5}$ and it is well known that there is almost surely an infinite coloured component if 
%and only if 
$\phi(p,q) > 0$.

% For the only if part, let $D_n$ (respectively $E_n$) be the event that there is a vertex $x$ (respectively, two vertices $x,y$) in $\partial B_n$ such that if $x$ is made red (respectively, $x$ and $y$ are made red) then there is a path in the graph from $B_{0.5}$ to $x$ consisting entirely of coloured vertices lying  inside $B_{n-0.5}$, except for $x$ itself. 

%These events  are measurable with respect to the $\sigma$-algebra generated by the Poisson process inside $B_{n-1/2}$ together with the colours (i.e. the variables $Y_i$ and $Z_i$) associated with the vertices inside $B_{n-1/2}$, and the locations only (but not the colours) of the vertices  in $\partial B_{n}$.  Let event $D_{n-1/2}$ be defined similarly.  If $A'_{n+1}$ occurs, then then at least one of events $D_n, D_{n-1/2}, E_n, E_{n-1/2}$ must occur.  But if $D_{n-1/2} \cup E_{n-1/2}$ occurs then the conditional probability of event $A_{n}$ occurring  is at least $p^3(1- \exp(-\lambda /16))$, this being $p^3$ times the probabiltiy that there is at least one Poisson point in the ball of radius $1/4$ centred at $(n-1/4)(x/|x|) $.

%If $(D_n \cup F_n) \setminus (D_{n-1/2} \cup E_n)$ occurs  then the conditional probability of $A_n$ is at least $p^2$, this being the probability that $x$ and $y$ are both red.  Combining these arguments yields 
%$$
%\theta_n(p,q) \geq p^3(1- \exp(-\lambda /16)) \phi_{n+1},
%$$
% If there is almost surely an infinite component then $\phi(p,q) >0$, and so $\theta(p,q) \geq  p^3(1- \exp(-\lambda /16)) \phi(p,q) >0$ 

 For the only if part, if there is almost surely an infinite component then $\phi(p,q) > 0$.   Given $n \geq 6$, we build up the Poisson process on the whole of $B_{n - 3}$. If there are any closed vertices that are not definitely correctly or incorrectly configured, we build up the process in the rest of their $1$-neighbourhood, and this determines whether they are green or uncoloured. If any more closed vertices occur they cannot be correctly configured as they will be joined to a closed vertex. Therefore we have built up the process everywhere in a region $R $ with $B_{n-3} \subset R \subset B_{n-2}$,
% boundary almost surely contained in $B_{n-2} \setminus B_{n-3}$
 and all uncoloured vertices at this stage will remain uncoloured.
 Let $V$ be the set of coloured vertices that are joined by a coloured path to a coloured vertex in $B_{0.5}$ at this stage.  

Next, we build out the process radially symmetrically from $B_{n-3}$ (apart from where the process has already been built up) until a vertex $v$ occurs that is connected to a vertex in $V$. Let $J$ be the
event that such a vertex $v$ occurs at distance $r$ between $n-3$ and
$n-1$ from the origin, so $J$ must occur for $A_n'$ to occur. 
We can find points $a_1, a_2, \ldots, a_9$ on the line $0 v$ extended away from the origin such that $a_1$ is 
$r + 0.3$ from the origin, $a_2$ is $r+0.6$ from the origin and so on. Surround $a_1,\ldots, a_9$ with circles $D_1, \ldots, D_9$ of
radius $0.05$ around them. If there is at least one red vertex in
each one of these little circles that is contained in $B_n$
when the process continues to the whole
of $B_n$, and $v$ is also red then $A_n$ occurs. Therefore if $J$ occurs then the
conditional probability of $A_n$ occuring is at least $\gamma$,
where
\[
    \gamma = p(1 - \exp(- 0.0025 \lambda p \pi))^9,
\]
as this is the probability of getting at least one red vertex in
each little circle and $v$ being red. Therefore $\theta_n (p,q) \geq \gamma P[J] \geq
\gamma \phi (p,q)$ for all $n \geq 6$, so $\theta (p,q) \geq \gamma
\phi (p,q) > 0$. \hfill{$\Box$} \vspace{.5cm}

%\begin{lemm}
% \[
% \frac{\partial \theta_n (p,q)}{\partial q} \geq \delta (p,q) \frac{\partial \theta_n (p,q)}{\partial p}
% \]

%where $\delta$ does not depend on $n$ and is strictly positive and continuous on $(0,1)^2$
%\end{lemm}

Our next lemma provides an analogue of the Margulis-Russo formula for
the enhanced continuum model. First, we need to introduce the notion of
pivotal vertices.

Given the configuration
$(\mathscr{P}_{\lambda},\underline{Y},\underline{Z})$ and
inserting a vertex at $x$ we say that $x$ is $1$-$pivotal$ $in$
$B_n$ if putting $Y_0 = 0$ means that $A_n$ occurs but putting $Y_0
= 1$ means it does not. Notice that $x$ can either complete a path
(but it cannot do via being enhanced), or it could make another
closed vertex correctly configured which in turn would complete a
path. We say that $x$ is $2$-$pivotal$ $in$ $B_n$ if inserting a
vertex at $x$ and putting $Z_0 = 0$ means $A_n$ occurs but putting
$Z_0 = 1$ means it does not. That is, $Y_0 > p$ and adding a closed
vertex $v$ at $x$ means $v$ is correctly configured and enhancing it
to a green vertex means $A_n$ occurs but otherwise it does not.

For $i = 1,2$ let $E_{n,i}(x)$ be the event that $x$ is $i$-pivotal in $B_n$,
and set $P_{n,i}(x,p,q) := P[E_{n,i}(x)]$.

\begin{lemm} \label{prop:Russo}
For all $n > 0.5$ and $p \in (0,1)$ and $q\in (0,1)$ it is the case that
\bea
\frac{\partial \theta_n(p,q)}{\partial p}=\int^{}_{B_n} \lambda P_{n,1}(x,p,q) \, \mathrm{d}x
%\eqno(1)
\label{eq1}
\eea
and
\bea
\frac{\partial \theta_n(p,q)}{\partial q} = \int^{}_{B_n} \lambda P_{n,2}(x,p,q) \, \mathrm{d}x.
%\eqno(2)
\label{eq2}
\eea
\end{lemm}
\noindent {\bf Proof.}
Let ${\cal F}$  be the $\sigma$-algebra generated
by the locations but not the colours of the vertices of $\Po_\lambda \cap B_n$.
Let $N_1$ be the number of $1$-pivotal vertices.
Define ${\cal F}$-measurable random variables,
$X_{p,q}$ and 
$Y_{p,q}$ as follows; $X_{p,q}$ is
the conditional probability that $A_n$ occurs,
 and $Y_{p,q}$ is the conditional expectation of $N_1$, 
 given the configuration
of $\Po_\lambda$. 
By the standard version of the Margulis-Russo formula
for an increasing event defined on a finite collection
of Bernoulli variables (Russo~(1981), Lemma 3),
$$
\lim_{ h \to 0}
h^{-1} (X_{p+ h,q} - X_{p,q})
= 
 Y_{p,q} , ~~ a.s.
$$
Let  $M$ denote the total number of vertices of $\Po_\lambda$ in $B_n$.
By the standard coupling of Bernoulli variables, and
Boole's inequality,
$|X_{p+ h,q} - X_{p,q}| \leq |h| M$ almost surely, 
and since $M$ is integrable we have by dominated convergence
that
\begin{eqnarray}
\frac{\partial \theta_n(p,q)}{
\partial p}
=
\lim_{ h \to 0}
E[ h^{-1} (X_{p+ \delta,q} - X_{p,q})]
= 
E[
 Y_{p,q} ] = E[N_1],
\label{Mar29a}
\end{eqnarray}
and by a standard application of the Palm theory of Poisson processes
(see e.g. Penrose~(2003)),
the right  hand side of (\ref{Mar29a})
equals the right hand side
of (\ref{eq1}). The proof of (\ref{eq2}) is similar. 
\hfill{$\Box$} \\

The key step in proving Theorem $1$ is given by the following result.
\begin{lemm}\label{intermediate}
There is a continuous function $\delta:(0,1)^2 \to (0,\infty)$ such that
for all $n > 100$, $x \in B_n$ and $(p,q) \in (0,1)^2$, we have
\begin{equation}
P_{n,2}(x,p,q) \geq \delta(p,q)P_{n,1}(x,p,q).
\label{100419a}
\end{equation}
% where $\delta$ is independent of $n$ and $x$ and is
% strictly positive and continuous on $(0,1)^2$.
\end{lemm}

Before proving this, we give a result saying that we can
assume there are only red vertices inside an annulus disk of fixed
size. For $x \in  \R^2$, and $0 \leq \alpha < \beta$,   let $C_{\alpha}(x)$ be the closed circle (i.e., disk) of radius $\alpha$
 around $x$, and let $A_{\alpha,\beta}(x)$ denote the annulus 
$C_{\beta}(x) \setminus  C_{\alpha}(x)$.
 Given $n$ and given $x \in B_n$, let 
$R_n(x,\alpha,\beta)$ be the event that all vertices in $A_{\alpha,\beta}(x) \cap B_n$
are red.

\begin{lemm}\label{intermed2}
Fix $\alpha >3$ and and $\beta > \alpha +3$.
There exists a continuous function $\delta_1:(0,1)^2 \to (0,\infty)$,
such that  
for all $(p,q) \in (0,1)^2$, all
  $n > \beta +3$ and all $x \in B_n$ with $|x|< \alpha -2 $ or 
$x > \beta +2$, we have
\[
P[E_{n,1}(x)  \cap R_n(x,\alpha, \beta) ] \geq \delta_1(p,q) P_{n,1}(x).
\]
\end{lemm}

\noindent {\bf Proof.}
 We shall consider a modified model, which is the same
as the enhanced model but with   enhancements suppressed
for all  those vertices lying in $A_{\alpha-1,\beta+1}(x)$.
Let $E'_{n,1}(x) $ be the event that $x$ is 1-pivotal in
the modified model.

Returning to the original model,
 we first create the Poisson process of intensity $\lambda$
in $B_n \setminus A_{\alpha-1,\beta +1}$,
and determine which of these vertices are red.
 Then we build up the Poisson process of intensity $\lambda$ inside
 $B_n \cap A_{\alpha-1,\beta +1} $
 and for any of these new vertices  
with more than $4$ neighbours, or
with at least one closed neighbour outside $A_{\alpha-1,\beta +1}$,  we decide
 whether they are red or closed. This decides whether or not they
 are coloured as these vertices cannot possibly become green because
 they are not correctly
 configured. We now can tell which of the closed vertices 
outside $A_{\alpha-1,\beta +1}$ are correctly configured, and we determine
which of these are green.

 This leaves a set $W$ of vertices inside $A_{\alpha-1,\beta +1}$ that have at most four neighbours. If we surround each vertex in $W$ by a circle of radius $0.5$ then we cannot have any point covered by more than $5$ of these circles as this means that there is a vertex in $W$ with at least $5$ neighbours. All
of these circles are contained in $C_{\beta+2}$,
 which has area $\pi(\beta+2)^2$. Therefore 
\[
|W| \leq \frac{5 \pi (\beta+2)^2}{0.5^2 \pi} =  20(\beta + 2)^2.
\]
For $x$ to have any possibility of being $1$-pivotal,
at this stage there must be a
 set $W'$ contained in $W$ such that if every vertex in $W'$ is coloured and every vertex in $W \setminus W'$ is uncoloured then $x$ becomes $1$-pivotal.
In this case, with probability at least
$[p(1-p)]^{20(\beta +2)^2}$ we have every vertex in 
$W'$ red and every vertex in $W \setminus W'$ closed, which would imply 
%corresponds to the
 event $E'_{n,1}(x)$
 occurring. Therefore
 $P[E'_{n,1}(x)] \geq [p(1-p)]^{20(\beta+2)^2} P[E_{n,1}(x)]$.

Now we note that the occurrence or otherwise of $E'_{n,1}(x)$
is unaffected by the addition or removal of closed vertices
in $A_{\alpha, \beta}(x)$. This is because the suppression
of enhancements in $A_{\alpha-1,\beta+1}$ means that 
these added or removed vertices cannot be enhanced themselves,
and moreover any vertices they cause to be correctly
or incorrectly configured also cannot be enhanced.

Consider creating the marked Poisson process in $B_n$,
with each Poisson point (vertex) $x_i$  marked with 
the pair $(Y_i,Z_i)$, in two stages. First, add all marked vertices
in $B_n \setminus A_{\alpha,\beta}(x)$, and just the red vertices
in $B_n \cap A_{\alpha,\beta}(x)$. Secondly, add the closed vertices
in $B_n \cap A_{\alpha,\beta}(x)$.
% (i.e., those  with $Y_i >p$).
 The vertices added at the second
stage have no bearing on the event $E'_{n,1}(x)$, so 
$E'_{n,1}(x)$ is independent of the event that no vertices at all are added in the second stage.
Hence, 
$$
P[E'_{n,1}(x) \cap R_n(x,\alpha,\beta)]
 \geq \exp(- (1-p) \lambda (\beta^2 - \alpha^2))
P[E'_{n,1}(x)], 
$$
with equality if $|x| \leq n- \beta$.

Finally, we use a similar argument to the initial argument in this proof.
Suppose
$E'_{n,1}(x) \cap R_n(x,\alpha,\beta)$ occurs. Then there exist at most 
$20(\beta +2)^2$ vertices in
 $A_{\beta,\beta +1}(x) \cup A_{\alpha -1,\alpha}(x)$
which are correctly configured for which the possibility of
enhancement has been suppressed. If we now allow these to be possibly enhanced,
there is a probability  of at least
$(1-q)^{20(\beta +2)^2}$ that none of them is enhanced, in which
case the set of coloured vertices is the same for the modified model
as for the un-modified model and therefore $E_{n,1}(x)$ occurs.  
Taking
 $$
\delta_1(p,q) = [p(1-p)(1-q)]^{20(\beta+2)^2} \exp(-(1-p) \lambda
 (\beta^2 - \alpha^2 )),
$$
 we are done. \hfill{$\Box$} \\

%So we take $\delta_1 := [p(1-p)(1-q)]^{20000}$ and we are done. \hfill{$\Box$} \\

\noindent {\bf Proof of Lemma \ref{intermediate}.}
As a start, we fix $p$ and $q$. 
 We also fix $n$ and  $x \in B_n$,  
and just write $P_{n,i}(x)$ for $P_{n,i}(x,p,q)$.
Define event $E_{n,1}(x)$
 as before,
so that $P_{n,1} (x) = P[E_{n,1}(x)]$. 
Also, write
$C_{r}$ for the disk $C_r(x)$.
% be the closed circle
%(i.e., disk) of radius $r$ around $x$.  
%so in previously used notation we have
%We alsow write $P'_{n,i}(x)$ for
%$P[E_1 \cap E_2]$.
For now we assume $30.5 < |x| < n-30.5$.
%$C_{30}$  does not intersect $B_{0.5}$ or $\partial B_n$.
 We create the Poisson process of intensity $\lambda$
 everywhere on $B_n$ except inside
$C_{30}$, and decide which of these vertices are red. 

%Then if there are any closed vertices that are not
%definitely correctly or incorrectly configured, we build up the
%process in the rest of their $1$-neighbourhood, and this determines
%whether they are coloured or uncoloured. If any more closed vertices occur
%they cannot be correctly configured as they will be joined to an
%uncoloured vertex. Therefore we have built up the process everywhere
%except for a region $R$ with boundary contained almost surely in $C_{30} \setminus
%C_{29}$, and all vertices so far will remain in the same state. Let
%$E_1$ be the event that there is no current coloured path from $B_{0.5}$
%to $\partial B_n$, so $E_1$ must occur for $x$ to be $1$-pivotal.

Now we create the process of only the red vertices 
 in $A_{25,30}$ (a Poisson process of intensity $p \lambda$
in this region).
Assuming there will be no closed vertices in $A_{25,30}$,
we then know which of the closed  vertices outside $C_{30}$ are correctly
configured, and we determine which of these are green.

Having done all this, let  $V$ denote the set of current  vertices 
in $B_n \setminus C_{25}$ that are connected to $B_{0.5}$ at this stage
(by connected we mean connected via a coloured path), and
let $T$ denote the set of current vertices in $B_n \setminus C_{25}$ that
 are connected to $\partial B_{n}$.

Let $N(V)$ be the  $1$-neighbourhood of $V$ 
and  let $N(T)$  be the  $1$-neighbourhood of $T$. 
 We build up the red process inwards (i.e., towards $x$ from the boundary
of $C_{25}$)
on $C_{25} \cap (N(V) \triangle N(T))$
 until a red vertex $y$
occurs (if such a vertex occurs). Set  $r= |y-x|$.
Suppose $y \in N(V)$ (if instead $y \in N(T)$
we would reverse the roles of $V$ and $T$ in the sequel).
Then if $T \cap C_{r+0.05} \neq \emptyset$ we say that 
event $F$ has occurred and we let $z$ denote an
arbitrarily chosen vertex
of  $T \cap C_{r+0.05} $.
 Otherwise, we  build up the red process inwards on
 $C_{r} \cap N(T) \setminus N(V)$ until a red vertex $z$
occurs (if such a vertex occurs). 
%Set $s = |z-x|$, so $ s \in ( r - 0.95, r+0.05)$.

Let $E_2$ be the event that (i)  such vertices $y$ and $z$ occur,
and (ii) the sets $V$ and $T$ are disjoint, and (iii)
 $|y-z| >1$, 
and (iv) there is no path from $y$ to $z$ through coloured
vertices in $B_n \setminus
C_{25}$ that are not in $V \cup T$.
If $E_{n,1}(x) \cap R_n(x,20,30)$ occurs, then $E_2$ must occur.

%Assume that the event $E_2$ occurs. We now build up the
%rest of the process on $C_{30} \setminus C_s $. Let $E_4$ be the event that
%there are no vertices of any kind in this region apart from the red
%vertices created previously. Then $P[E_4 | E_2] \geq \exp
%(-\lambda \pi 31^2) =: \delta_2$.

\begin{figure}[htbp]
\includegraphics[angle = 270, width = 14cm]{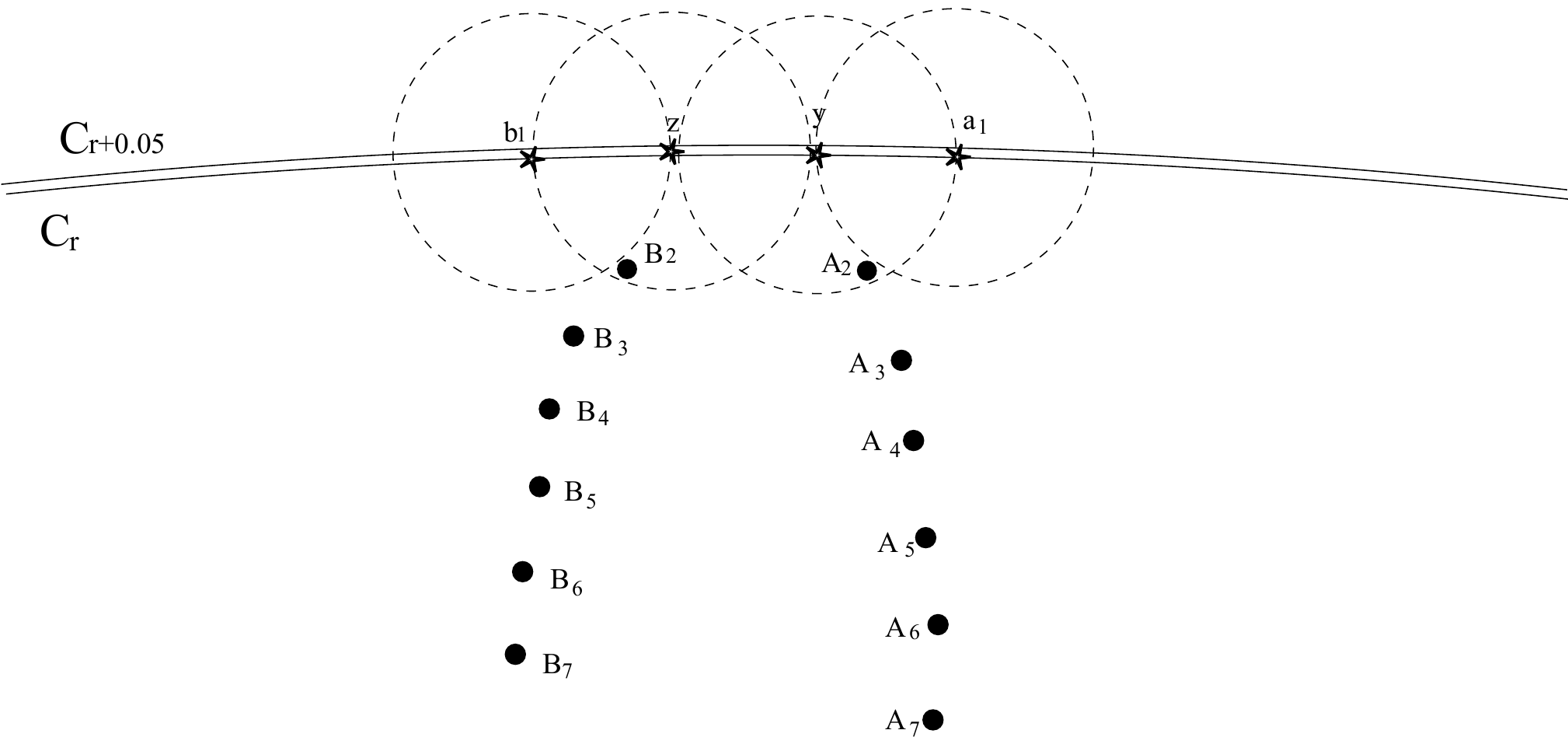}
\caption{Our convention in the diagrams is to indicate points with
lower case letters, and areas with upper case letters. The
dashed circles are of radius 1. Here the event $F$ occurs. 
} \label{fig:fig1}
\end{figure}

Now suppose $E_2 \cap F$ has occurred.  
%We first assume that $F$ occurs.
Let $a_1$ be the point (again we use `point'
to refer to a point in $\mathbb{R}^2$) which is at distance
$r$ from $x$ and distance $1$ from $y$ on the opposite side of the
line $xy$ to the side $z$ is on (see Figure $\ref{fig:fig1}$). 
Let $a_2$ be the point lying inside $C_r$ at distance 
$1.01$ from $a_1$ and $0.99$ from $y$.
 Let $A_2$ be the circle of radius $0.005$ around $a_2$.
Any red vertex in this circle will be connected to the red vertex
$y$ (and therefore to a path to $B_{0.5}$) but cannot be connected
to any coloured path to $\partial B_n$ as $a_1$ is the nearest place for
such a vertex to be, given $E_2$ occurs. Similarly let 
$b_1$ be the point lying at distance $1$ from $z$ and distance $r$ from
$x$, on the opposite side of $xz$ to $y$.  Then let
 $b_2$ be the point at distance $1.01$ from $b_1$ and $0.99$ from
$z$, and let $B_2$  be the circle of radius $0.005$ around $b_2$. Any red
vertex in $B_2$ will be connected to $z$ (and therefore a path to
$\partial B_n$), but not a path to $B_{0.5}$. Also, any vertex in
$A_2$ will be at least $1.1$ away from any vertex in $B_2$. 

Now let $l$ be the line through $x$ such that $a_2$ and $b_2$ are on
different sides of the line and at equal distance from the line. We
can pick points $a_3,a_4,\ldots,a_{30}$ such that $a_i$ is within $0.9$
of $a_{i + 1}$ for $2 \leq i \leq 29$, $a_{30}$ and $a_{29}$ are
both within $0.9$ of $x$ but none of the other $a_i$'s are within $1.1$ of $x$, and none of the ${a_i :i \geq 3}$ are within
$1$ of $C_r$ or within $0.5$ of $l$ or within $0.01$ of another
$a_j$. Do the same on the other side of $l$ with $b_3,b_4, \ldots , b_{30}$.
Now consider circles $A_i$ and $B_i$ of radius $0.005$ around them.
Let $I$ be the event that there is at exactly one red vertex in each
of these circles, and also the circles $A_2$ and $B_2$, and there
are no more new vertices anywhere else in $C_{25}$, and no closed
vertices in $C_{30}\setminus C_{25}$. The probability that $I$
occurs, given $E_2 \cap F$, is at least
\[
 \delta_2 := (1 - \exp(-0.005^2 \pi \lambda p))^ {60} \exp(-900 \pi \lambda).
\]
 If the events $E_2$,  $F$, $I$ occur and $Y_0
> p$ then $x$ is $2$-pivotal.

\begin{figure}[htbp]
\includegraphics[angle = 270, width = 12cm]{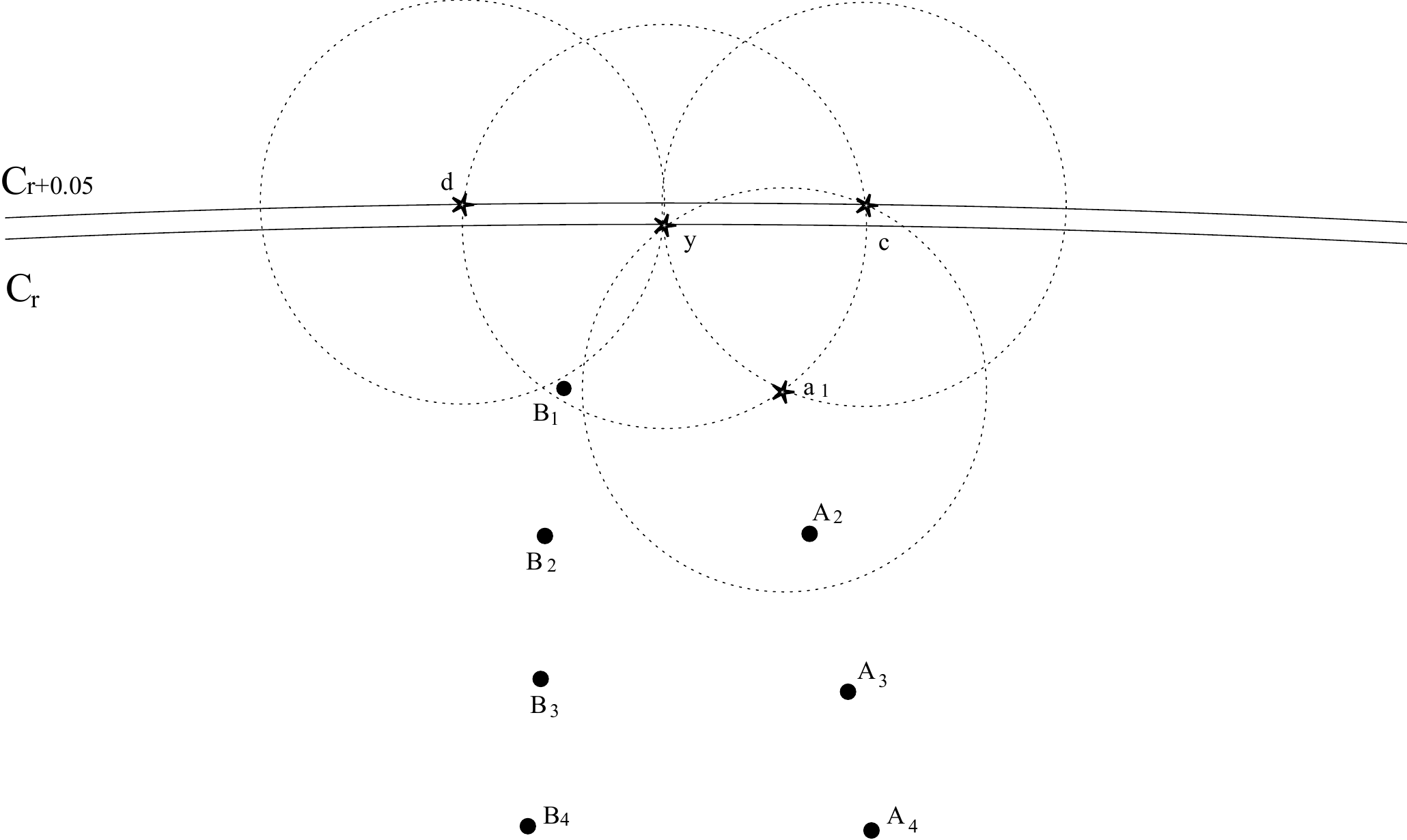}
\caption{The case where $F$ does not occur. Here $a_1$ is the `worst possible'
location for $z$}
\label{fig:points1}
\end{figure}

Now we consider the case where $E_2$ occurs but $F$ does not, so $z$ is inside
$C_r$ and is connected to a vertex $z_1$ in $T$ that must be outside
$C_{r+0.05}$ as $z$ is the only vertex in $T$ inside $C_{r+0.05}$ (see Figure $\ref{fig:points1}$).

 Let $c$ be the point at distance $1$ from $y$ and $r+0.05$ from $x$, on
the same side of the line $xy$ as $z$ (assume without loss of
generality this is to the right of $y$). This is the closest $z_1$
can be. Let $a_1$ be the point inside $C_r$ at distance $1$ from $y$
and $1$ from $c$, so this is the furthest left that $z$ can be. Let
$d$ be the point at distance $r+0.05$ from $x$ and $1$ from $y$, on
the other side of $y$ to $c$. Then consider the point $b_1$ inside $C_r$ at
distance $1.01$ from $d$ and $0.99$ from $y$, and the small circle
$B_1$ of radius $0.005$ around $b_1$. Then any vertex in $B_1$ is distant at
least $1.01$ from $a_1$, and therefore from $z$, as $z$ cannot be any
nearer than $a_1$. Also any vertex in $B_1$ will be at least $1.005$
from any other vertices in $T$, as $d$ is the nearest place such a
point can be. As before we can then have points $a_2,\ldots,a_{30}$
and $b_2,\ldots,b_{30}$ with small circles around them such that having one red vertex in each of these vertices ensures that $x$ is $2$-pivotal. The
probability of getting at least $1$ red vertex in each of these
circles, a red vertex in $B_1$ and no other new vertices in $C_{25}$,
and no closed vertices in $C_{30} \setminus C_{25}$, is
at least $\delta_2$.

So by Lemma \ref{intermed2},
 the probability that $x$ is $2$-pivotal satisfies
\beaa
    P_{n,2}(x)&  \geq &   \delta_2 P[E_2  \cap F] +
      \delta_2 P[E_2 \cap F^c] 
\\ 
& \geq &  
%\delta_2
%    \delta_3 P[E_3]
  %\geq
 \delta_2  P[ E_{n,1}(x) \cap R_n(x,20,30)]
\\ & \geq & \delta_1 \delta_2  P_{n,1}(x)
.
\eeaa
 This proves the claim (\ref{100419a}) for the case with $30.5 < |x| < n-30.5$.

Now suppose $|x| \leq 30.5$. Then
we create the Poisson process  in $B_n \setminus C_{40}$,
 and decide which of these vertices are red.  Then  we create  the
red process in  $A_{39,40}(x)$, and determine which
vertices in $B_n \setminus C_{40}$ are green, assuming
there are no closed vertices in $A_{39,40}(x)$.
 We then build up the red process in $C_{39}$ inwards towards $x$
until a vertex $y$ occurs in the process which is connected to
$\partial B_n$. Let $H_1$ be the event that such a vertex $y$
appears at distance $r$ between $38$ and $39$ from $x$, so $H_1$ must
occur for $E_{n,1}(x) \cap R_n(x,20,40) $ to occur.

If $x$ is inside $B_{0.5}$ we can choose points $a_0$ and $a_1$ such
that they are both outside $B_{0.5}$, at distance between $0.8$ and
$0.9$ from $x$ and at distance between $0.1$ and $0.2$ from each
other. We can then choose $b_0$ and $b_1$ such that they are both
within $0.9$ of $x$, further than $1.5$ from $a_0$ and $a_1$ and
between $0.1$ and $0.2$ from each other. We can then choose points
$a_2, a_3, \ldots, a_{100}$ such that $a_i$ is within $0.9$ of $a_{i+1}$
for $1 \leq i \leq 99$, $a_{100}$ is within $0.9$ of $y$, no two
$a_i$ are within $0.1$ of each other, and no $a_i$ is within $1.1$
of $x$, $b_0$ or $b_1$, or inside $B_{0.5}$ for $i \geq 2$. Then
consider little circles $A_i$ and $B_i$ of radius $0.05$ around
these points. If there is at least one red vertex in each of these
circles and no vertices anywhere else in $C_r$ then $x$ is
$2$-pivotal. If $x$ is outside $B_{0.5}$ we choose points in a
similar way but make sure $b_1$ connects with a path to $B_{0.5}$,
using little circles $B_2, B_3,\ldots,B_{50}$ which are again of radius
$0.05$ and are at least $1.1$ from the $A_i$. Therefore, 
setting
\[
\delta_3 := (1-\exp(-0.05^2 \pi \lambda p))^{152}\exp(-1600 \pi
\lambda)
\]
and using
Lemma \ref{intermed2}, we have for some strictly positive continous
$\delta_4(p,q)$ that
%probability that $x$ is $2$-pivotal is at least 
$$
P_{n,2}(x) \geq \delta_3 P[H_1] \geq \delta_3 P[E_{n,1}(x)
  \cap
R_n(x,20,40)] \geq  \delta_3 \delta_4 P_{n,1}(x).
$$

Now suppose $|x| \geq n - 30.5$.
 In this case the proof is similar. Again, create
 the Poisson process in $B_n \setminus C_{40}$.
 Then create  the red process in $A_{39,40}(x)$ and
determine which vertices in $B_n \setminus C_{40}$ are green,
assuming there are no closed vertices in $A_{39,40}(x)$.
Then  build the red process in $C_{39}  \cap  B_{n-0.5}$
  inwards towards $x$
%  (but excluding $\partial B_n \cap C_{40}$) 
until a vertex $y$ occurs that is connected to a
path of coloured vertices to $B_{0.5}$ but not to $\partial B_n$. Let
$H_2$ be the event that such a vertex $y$ occurs at distance $r$
between $38$ and $39$ from $x$, and that there is no current coloured
path from $B_{0.5}$ to $\partial B_n$, so $H_2$ has to occur for 
$E_{n,1}(x) \cap R_n(x,20,40)$ to occur.
 Given this vertex $y$
 we can find circles
$A_1, A_2, \ldots ,A_{100}$ and $B_1, B_2, \ldots,B_{50}$ of radius $0.05$
as before such that having a red vertex in each of these little
circles but no other vertices in $C_r$ or $\partial B_n \cap C_{40}$
ensures $x$ is $2$-pivotal. Therefore in this case
\[
P_{n,2}(x) \geq  \delta_3 P[H_2 ] \geq \delta_3 P[E_{n,1}(x) \cap
R_n(x,20,40)] \geq \delta_3  \delta_4
P_{n,1}(x).
\]

Take $\delta(p,q) := \delta_1 \delta_2 \delta_3 \delta_4 $. By its
construction $\delta$ is strictly positive and continuous in $p$ and $q$,
completing the proof of the lemma. \hfill{$\Box$} \vspace{.5cm}

The following proposition follows immediately by combining 
Lemma~\ref{prop:Russo} and Lemma~\ref{intermediate}.
\begin{prop} \label{prop:final}
There is a continuous function $\delta:(0,1)^2 \to (0,\infty)$
such that
\[
\frac{\partial \theta_n (p,q)}{\partial q} \geq \delta (p,q) \frac{\partial \theta_n (p,q)}{\partial p}
\]
for all $n \geq 100$ and $(p,q) \in (0,1)^2$.
%it is the case that
%where $\delta$ does not depend on $n$ and is strictly positive and continuous on $(0,1)^2$.
\end{prop}

\noindent
\textbf{Proof of Theorem~\ref{bool}.}
Set $p^* = p_c^{\rm site} $
and $q^* = (1/8)(p^*)^2$. 
%and $\epsilon \in \min(0,
%For any $q$,
Then using Proposition~\ref{prop:final} and looking at a small
 box around $(p^*,q^*)$, we can find 
$\epsilon \in (0,\min(p^*/2,1-p^*))$ 
 and $ \kappa \in (0,q^*)$ 
such that for all $n > 100$ we have
\[
%\theta_n(p_c^{\rm site} + \epsilon,q^* - \kappa) \leq
% \theta_n(p_c^{\rm site} - \epsilon,q^* + \kappa).
\theta_n(p^* + \epsilon,q^* - \kappa) \leq
 \theta_n(p^*  - \epsilon,q^* + \kappa).
\]
Taking the limit inferior as $n \rightarrow \infty$,
since
 $\theta$ is monotone in $q$ we get
\[
0 < \theta(p^* + \epsilon,0 ) \leq
\theta(p^*  + \epsilon,q^*- \kappa)
 \leq \theta(p^* - \epsilon,q^* + \kappa).
\]
Now set $p = p^* - \epsilon $. Then   $q^* + \kappa \leq p^2$,
% we have can now find $p$ strictly
% less than $p_c^{\rm site}$ such 
so that $\theta(p,p^2)>0$, and
 by Proposition
\ref{prop:enh}, the enhanced model with parameters $(p,p^2)$ percolates,
i.e. has an infinite coloured component, almost surely.

We finish the proof with a coupling argument along the lines of
Grimmett and Stacey~(1998). Let E be the set of edges and $V$ be the
set of vertices of $\CC$ (the infinite component). Let $(X_e: e \in
E)$ and $(Z_v: v \in V)$ be collections of independent Bernoulli
random variables with mean $p$. From these we construct a new
collection $(Y_v: v \in V)$ which constitutes a site percolation
process on $\CC$. Let
 $e_0, e_1, ...$ be an enumeration of the edges of $\CC$ and
 $v_0, v_1, ...$ an enumeration of the vertices. Suppose
at some point we have defined $(Y_v: v \in W)$ for some subset $W$
of $V$. Let $\cal Y$ be the set of vertices not in $W$ which are
adjacent to some currently active vertex (i.e. a vertex $u \in W$
with $Y_u = 1$). If $\cal Y = \emptyset$ then let $y$ be the first
vertex not in $W$ and set $Y_y = Z_y$ and add $y$ to $W$. If $\cal Y \neq \emptyset$,
we let $y$ be the first vertex in $\cal Y$ and let $y'$ be the first
currently active vertex adjacent to it, then set $Y_y = X_{yy'}$ and add $y$ to $W$.
Repeating this process builds up the entire red site percolation process,
if it does not percolate, or a percolating subset of the red site percolation
process if it does percolate.
%(i.e., has an infinite component). 
In the latter case, the bond process $\{X_e\}$ also percolates.

Now suppose the red site process does not percolate.
For any correctly configured vertex $x_1$ with $x_2$ up to $x_5$ as
before, $x_1$ itself is not red. Therefore at most one edge to $x_1$
has been examined, so we can can find a first unexamined edge (in
the enumeration) to $x_2$ or $x_3$, and then to $x_4$ or $x_5$. We
then declare $x_1$ to be green only if both of these edges are open,
which happens with probability $p^2$.
This completes the enhanced site process with
 $q = p^2$ and every component
in this is contained in a component for the bond process $\{X_e\}$.

Therefore, since the enhanced $(p,p^2)$ site process percolates almost
surely,
so does the bond process, 
so $p_c^{\rm bond} \leq p < p_c^{\rm site}$. 
\hfill{$\Box$}
\vspace{.5cm}

%\section{Bibliography}

\end{document}